\documentclass{article}

\usepackage{geometry}
\usepackage{amsmath, amsfonts, bm, amssymb, enumitem, tabularx}
\usepackage{amsfonts, graphicx, grffile,fullpage, float, latexsym, mathtools, bbm}
\usepackage{hyperref}
\hypersetup{
    colorlinks=true, 
    citecolor=black,
     linkcolor=blue,
     urlcolor=blue,
    linktoc=all,     
    linkcolor=black,  
}
\usepackage{cleveref}
\usepackage{multicol}
\usepackage{tikz}
\usepackage{subcaption}

\usepackage[sort, numbers, compress]{natbib}

\usepackage{Shorthands}

\begin{document}

\title{Augmented Lagrangian Methods as Layered Control Architectures}
\author{Anusha Srikanthan{$^\dagger$}, Vijay Kumar{$^\dagger$}, Nikolai Matni{$^\dagger$\thanks{This research is in part supported by NSF awards CPS-2038873, SLES-2331880, and NSF CAREER award ECCS-2045834.}}, \\
{$^\dagger$}University of Pennsylvania\\
\tt\small \{sanusha, kumar, nmatni\}@seas.upenn.edu%
}

\maketitle

\begin{abstract}
    For optimal control problems that involve planning and following a trajectory, two degree of freedom (2DOF) controllers are a ubiquitously used control architecture that decomposes the problem into a trajectory generation layer and a feedback control layer.  However, despite the broad use and practical success of this layered control architecture, it remains a design choice that must be imposed \emph{a priori} on the control policy.  To address this gap, this paper seeks to initiate a principled study of the design of layered control architectures, with an initial focus on the 2DOF controller.  We show that applying the Alternating Direction Method of Multipliers (ADMM) algorithm to solve a strategically rewritten optimal control problem results in solutions that are naturally layered, and composed of a trajectory generation layer and a feedback control layer.  Furthermore, these layers are coupled via Lagrange multipliers that ensure dynamic feasibility of the planned trajectory.  We instantiate this framework in the context of deterministic and stochastic linear optimal control problems, and show how our approach automatically yields a feedforward/feedback-based control policy that exactly solves the original problem.  We then show that the simplicity of the resulting controller structure suggests natural heuristic algorithms for approximately solving nonlinear optimal control problems.  We empirically demonstrate improved performance of these layered nonlinear optimal controllers as compared to iLQR, and highlight their flexibility by incorporating both convex and nonconvex constraints. 
\end{abstract}

\section{Introduction}

Optimal control has proven to be a key approach to solving problems across a wide range of fields, including economics, robotics, and communication systems. However, despite their significance, solving optimal control problems can be challenging due to nonlinear dynamics, high-dimensional state and control spaces, uncertainty, noise, and constraints.  For optimal control problems that involve planning and following a trajectory, a ubiquitous \emph{layered control architecture}~\cite[Ch. 15]{aastrom2021feedback} commonly referred to as a \emph{two degree of freedom (2DOF) controller}~\cite{murray2009optimization} has emerged as the standard solution approach.  This control architecture decomposes the problem into a trajectory generation layer, which generates the nominal trajectory the system should follow, and a feedback control layer, which corrects for errors between the actual system evolution and the planned trajectory.  Indeed, this control architecture can be observed across linear control (feedforward/feedback control), robust model predictive control, and nonlinear control, and has led to significant practical impact across a wide variety of fields including robotics, power systems, communication networks, and biology.  

We elaborate more on these different settings below, but highlight here that despite the ubiquity and practical success of the layered approach, this control structure does not emerge naturally from solving an optimal control problem, but rather must be imposed \emph{a priori} on the control policy.  To address this gap, we seek to initiate a principled study of the design of layered control architectures, with an initial focus on the 2DOF design pattern.  Our starting point is the observation that Augmented Lagrangian-based optimization algorithms applied to optimal control problems can be naturally interpreted as two degree of freedom layered control architectures.  We instantiate this observation in the context of linear and nonlinear optimal control problems, and show perhaps surprisingly that solutions obtained using the Alternating Direction Method of Multipliers (ADMM) algorithm~\cite{boyd2011distributed} to solve the original optimal control problem are naturally layered and composed of a trajectory generation layer and feedback control layer.  In contrast to ad-hoc designs however, these two layers are coupled via Lagrange multipliers which ensure consistency between the planned trajectory and the tracking ability of the closed-loop feedback control layer.

\emph{Contributions:}
This paper seeks to initiate the study of layered control architectures (LCAs) through the lens of optimization algorithms.  Our specific contributions are:
\begin{enumerate}[leftmargin=*]
    \item We show that strategically applying the ADMM algorithm to solve an optimal control problem results in a natural 2DOF layered control architecture composed of a trajectory generation layer and a feedback control layer.  Importantly, the two layers are coupled via Lagrange multipliers that ensure dynamic feasibility of the planned trajectory.
    \item In the case of linear optimal control problems with convex but otherwise arbitrary cost functions, we show how this approach \emph{automatically yields} a feedforward/feedback controller that exactly solves the original problem.  We also show how this perspective allows us to seamlessly incorporate stochastic process noise into the problem.
    \item In the case of nonlinear optimal control problems, we exploit the structural simplicity of the resulting controller to propose a heuristic algorithm for \emph{constrained nonlinear optimal control} that uses iLQR~\cite{todorov2005generalized} as a sub-routine.  Although not the main focus of the paper, we emphasize the exciting possibilities that this novel perspective raises for nonlinear control design.
    \item We provide empirical evaluations that demonstrate the benefits of layered control strategies in the context of nonlinear optimal control.
\end{enumerate}

\emph{Related work---2DOF and Layered Control Architectures:} 
In linear control systems, 2DOF controllers decompose the control input into a feedforward term, which drives the system to the desired trajectory, and a feedback term, which compensates for errors~\cite{murray2009optimization}.  Analogous design patterns are observed in robust model predictive control (MPC)~\cite{bemporad2007robust}.  For example, tube-based MPC approaches~\cite{borrelli2017predictive} broadly apply a control input of the form $u = K(x-x_d) + u_d$, where $(x_d,u_d)$ are nominal state and control inputs computed by solving an optimization problem online, and $K(x-x_d)$ is a feedback term compensating for errors between the actual system state $x$ and the reference state $x_d$. 

For nonlinear systems, trajectory generation and feedback control are typically decoupled, although approaches exist that do not explicitly make this separation, e.g., iLQR~\cite{todorov2005generalized}.  A typical design pattern consists of generating a(n) (approximately) dynamically feasible and safe reference trajectory, e.g., by exploiting differential flatness or a reduced order model, and then applying locally stabilizing feedback control to ensure trajectory tracking, e.g., via linearization or control Lyapunov functions.  Recent efforts from the robotics community show how to obtain ``full-stack'' safety/stability/performance guarantees for such layered architectures, see for example~\cite{rosolia2020multi,csomay2022multi,rosolia2022unified,garg2021multi}, by appropriately constraining planned trajectories to account for feedback control tracking error.  In addition, work from formal methods solving discrete planning problems over continuous dynamics can be viewed as a layered approach, to solving a complementary planning and control problem~\cite{lindemann2018control,dimitrova2014deductive,raman2014model,sadraddini2015robust,wolff2014optimization,wolff2012robust,fan2020fast,kress2009temporal}.  This body of work is exciting, as it treats layered control architectures as an object of study, and provides formal guarantees of correctness.  We emphasize however that these papers impose the layered architecture \emph{a priori}, and as such, do not address the question of how such layered architectures can be derived from first principles.

\emph{Related work---Theory of Layered Architectures:}
Originally motivated by communication networks~\cite{chiang2007layering}, the Layering as Optimization Decomposition (LAO)~\cite{palomar2006tutorial} perspective has emerged as a promising quantitative theory of layered architectures. At a high-level, the LAO framework argues that layered architectures can be viewed as arising from a vertical decomposition of an optimization problem, wherein redundant variables are introduced across layers, and coordination enforced via Lagrange multipliers.  LAO has been successfully applied to both communication~\cite{chiang2007layering} and power systems~\cite{zhao2014design,cai2017distributed}, resulting in exciting breakthroughs in both fields.   This work however focused on the solution of static optimization problems, e.g., Network Utility Maximization or Optimal Power Flow problems.  To the best of our knowledge, the first extension of these ideas to optimal control problems can be found in~\cite{matni2016theory}, where a LAO inspired relaxation is applied to a distributed linear optimal control problem in order to obtain a layered control architecture that approximately solves the original optimal control problem.  A main contribution of this work was the \emph{derivation} of a dynamics-aware trajectory planning layer, wherein the trajectory planning problem is augmented with a tracking penalty that characterizes the feedback control layer's ability to follow a given trajectory. We then extend this approach to nonlinear systems in~\cite{srikanthan2023data}, where we propose a data-driven approach to approximating the aforementioned tracking penalty for a fixed feedback controller. We note however that in both~\cite{matni2016theory, srikanthan2023data}, layered architectures are only obtained by considering suitable \emph{relaxations} of the original optimal control problem.  In contrast, in this work we show how optimization algorithms used to directly solve the original problem can be interpreted as layered control architectures themselves.

\emph{Paper organization:} 
We show how ADMM applied to an optimal control problem results in a layered control architecture in Section~\ref{sec:prelims}. In Section~\ref{sec:method}, we instantiate our layered control architecture in the context of a deterministic and stochastic linear system. In Section~\ref{sec:non-lin}, we discuss extensions to the nonlinear setting.  We evaluate our proposed approach in several numerical examples in Section~\ref{sec:sim}, and end with conclusions and future work in Section~\ref{sec:conclusion}.

\section{Problem Formulation}\label{sec:prelims}

We consider the discrete-time finite-horizon optimal control problem (OCP) with initial condition $x_0 = \xi$
\begin{equation}\label{prob:ocp}
\begin{array}{rl}
     \mathrm{minimize}_{\mathbf{x, u}} &\mathcal{C}_x(\mathbf{x}) + \mathcal{C}_u(\mathbf{u}) \\
     \text{s.t.} &  x_{t+1} = f(x_t, u_t)\\
     & x_t \in \mathcal{R}, \ t=0, \dots, N-1\\
     & x_0 = \xi
\end{array}
\end{equation}
where $x_t \in \mathbb{R}^n$ is the state, $u_t \in \mathbb{R}^m$ is the control input, $\mathbf{x}:=(x_0,x_1,\dots,x_{N}) \in \mathbb{R}^{n \times N+1}$ is the state trajectory, $\mathbf{u}:=(u_0,u_1,\dots,u_{N-1}) \in \mathbb{R}^{m \times N}$ is the control input trajectory, $C_x(\mathbf{x})$ is the state cost, $C_u(\mathbf{u})$ is the input cost, $\mathcal{R} \subseteq \mathbb{R}^n$ is a state constraint set, and $f$ is the nonlinear dynamics function. 

While many approaches to solving OCP~\eqref{prob:ocp} exist~\cite{murray2009optimization}, our goal is to define a solution strategy which systematically generates 2DOF layered control architecture. Towards that end, we consider the equivalent OCP
\begin{equation}\label{eq:redundant-ref-var}
\begin{array}{rl}
    \mathrm{minimize}_{\mathbf{r, x, u}} & \mathcal{C}_x (\mathbf{r}) + \mathcal{C}_u (\mathbf{u}) \\
     \text{s.t.} & r_t \in \mathcal R\\
     &x_{t+1} = f(x_t, u_t), \ t=0, \dots, N-1 \\
     & x_0 = \xi \\
        &  \mathbf{r} = \mathbf{x} 
\end{array}
\end{equation}
obtained from~\eqref{prob:ocp} through the introduction of a redundant ``reference variable'' $\vec r = (r_0,r_1,\dots,r_{N})$ constrained to satisfy $\vec r = \vec x$.
We now show how solving~\eqref{eq:redundant-ref-var} using ADMM naturally yields solutions with a 2DOF layered control architecture.

\subsection{Alternating direction method of multipliers}

The following is adapted from~\cite{boyd2011distributed}.
Consider the optimization problem
\begin{equation}\label{eq:admm_orig}
    \begin{aligned}
    \mathrm{minimize} &\quad g(r) + h(z) \\
    \mathrm{subject\ to} &\quad Ar + Bz = c,
    \end{aligned}
\end{equation}
over the decision variables $r$ and $z$, with convex functions $f$ and  $g$. Define the scaled-form augmented Lagrangian of optimization problem~\eqref{eq:admm_orig} as:
\begin{equation}
    L_{\rho}(r, z, v) = g(r) + h(z) + \frac{\rho}{2} \Vert Ar + Bz - c  + v\Vert_2^2 - \frac{\rho}{2}\|v\|^2_2
\end{equation}
where $(r, z)$ are the primal variables, $v$ is the (scaled) dual variable associated with the equality constraint, and $\rho>0$ is an algorithm parameter. The constrained optimization problem~\eqref{eq:admm_orig} is solved by alternatively minimizing the scaled-form augmented Lagrangian over the primal variables $r$ and $z$, and updating the (scaled) dual variable $v$:
\begin{equation}\label{alg:admm}
    \begin{array}{rcl}
        r^{k+1} &\coloneqq& \argmin_x L_{\rho} (r, z^k, v^k) \\
        z^{k+1} &\coloneqq& \argmin_z L_{\rho} (r^{k+1}, z, v^k) \\
        v^{k+1} &\coloneqq& v^k + (Ar^{k+1} + Bz^{k+1} - c)
   \end{array}
\end{equation}

Next, we describe the convergence properties of ADMM. Suppose optimization problem~\eqref{eq:admm_orig} satisfies the two assumptions stated below, then the following theorem holds.

\begin{assumption}\label{a1}
    The (extended-real-valued) functions $f: \mathbb{R}^n \longrightarrow \mathbb{R} \cup \{+\infty\}$, and $g: \mathbb{R}^m \longrightarrow \mathbb{R} \cup \{+\infty\}$ are closed, proper, and convex.
\end{assumption}

\begin{assumption}\label{a2}
    The standard Lagrangian for problem \eqref{eq:admm_orig} has a saddle point. 
\end{assumption}

\begin{theorem}[\S 3.2.1 in \cite{boyd2011distributed}]\label{thm:admm}
Let $p^\star$ denote the optimal value of optimization problem~\eqref{eq:admm_orig}. Under Assumptions $1$ and $2$, the ADMM iterates satisfy the following:

\begin{itemize}
    \item \textit{Residual convergence}: $Ar^{k} + Bz^{k} - c \to 0$ as $k \to\infty$, i.e., the iterates approach feasibility.
    \item \textit{Objective convergence}: $f(r^k) + g(z^k) = p^*$ as $k \to \infty$, i.e., the objective function of the iterates approaches the optimal value.
    \item \textit{Dual variable convergence}: $v^k \to v^\star$ as $k \to \infty$, where $v^\star$ is a dual optimal point. 
\end{itemize}
\end{theorem} 

Finally, we note that ADMM has been widely applied to solve nonconvex optimization problems.  Rapid convergence to local optima has been observed empirically in a variety of settings, and can be guaranteed under certain assumptions~\cite{liu2019linearized}.

\subsection{ADMM yields 2DOF layered control architectures}

The ADMM iterates~\eqref{alg:admm}, when instantiated on OCP~\eqref{eq:redundant-ref-var}, become
\begin{subequations}
    \begin{align}
        \mathbf{r}^{k+1} &\coloneqq \argmin_{\mathbf{r}} \
         C_x(\mathbf{r}) + \frac{\rho}{2}\Vert \mathbf{x}^k - \mathbf{r} + \mathbf{v}^k \Vert_2^2 \nonumber \\
        & \quad \text{s.t. }\mathbf{r} \in \mathcal{R}^N \label{eq:ref-layer} \\
        \left(\mathbf{x}^{k+1}, \mathbf{u}^{k+1}\right) &\coloneqq \argmin_{\mathbf{x, u}} \ \frac{\rho}{2}\Vert \mathbf{x} - \mathbf{r}^{k+1} + \mathbf{v}^k \Vert_2^2 + C_u(\mathbf{u})  \nonumber \\
         & \quad\text{s.t. } x_{t+1} = f(x_t, u_t), \ t=0, \dots, N-1 \label{eq:track-layer} \\
         & \quad\quad \ \, x_0 = \xi \nonumber \\
        \mathbf{v}^{k+1} &\coloneqq \mathbf{v}^k + \mathbf{x}^{k+1} - \mathbf{r}^{k+1} \label{eq:dual-update}
\end{align}
\label{eq:admm_ocp}
\end{subequations}
where $\mathcal{R}^N:=\mathcal R \times \cdots \times \mathcal R$ is the Cartesian product of the  constraint set $\mathcal{R}$ over the time horizon $N$.


We describe how the ADMM iterate updates~\eqref{eq:admm_ocp} can be interpreted as a layered control architecture:
\begin{enumerate}[label=\alph*)]
    \item \textbf{Trajectory generation layer~\eqref{eq:ref-layer}:} The $\vec r$-update step~\eqref{eq:ref-layer} is naturally interpreted as a trajectory generation layer, wherein an updated reference trajectory $\mathbf{r}$ is obtained by optimizing the utility cost $C_x(\vec r)$ subject to state constraints $\vec r \in \mathcal R^N$.  While the reference trajectory is not explicitly constrained to be dynamically feasible a trust-region-like penalty $\frac{\rho}{2}\Vert \mathbf{x}^k - \mathbf{r} + \mathbf{v}^k \Vert_2^2$ arising from the augmented Lagrangian regularizes the reference trajectory to be approximately consistent with the current (dynamically feasible) state trajectory $\vec x^k$.
    \item \textbf{Feedback control layer~\eqref{eq:track-layer}:} We immediately recognize the $(\vec x, \vec u)$-update step~\eqref{eq:track-layer} as a reference tracking optimal control problem, with reference given by $\vec r^{k+1} - \vec v^k$.  Depending on the problem setting, exact or approximate optimal feedback controllers can be obtained to this update step.
    \item \textbf{Dual update~\eqref{eq:dual-update}:} Finally, the dual variables $\vec v$ are updated according to equation~\eqref{eq:dual-update}.  We notice that the dual variable $\vec v$ can be seen as a protocol between the trajectory generation and feedback control layers that ensure that the planned reference trajectories converge to dynamically feasible behaviors (and vice versa).
\end{enumerate}

In the next sections, we instantiate this framework in the context of linear and nonlinear optimal control problems.  For linear optimal control problems with convex costs, we show that the solution produced by the updates in~\eqref{eq:admm_ocp} is a 2DOF optimal controller with a trajectory generator along with feedforward and feedback control terms.  For nonlinear optimal control problems, we show a natural separation between planning and control that isolates challenging lower-layer nonlinear feedback control from higher-layer trajectory generation and planning. 

\section{Layered Control Architectures for Linear Systems}\label{sec:method}

In this section, we instantiate the ADMM updates~\eqref{eq:admm_ocp} in deterministic and stochastic linear OCPs, and show convergence to the optimal solution when the cost function is convex.

\subsection{Deterministic linear system}\label{sec:no-noise-lin}

We consider the deterministic linear dynamics $x_{t+1} = A_t x_t + B_t u_t$, and quadratic control cost $C_u(\mathbf{u}) = \sum_{t=0}^{N-1} u_t^T R_t u_t$, for $R_t$ positive definite matrices. The ADMM updates~\eqref{eq:admm_ocp} then become
\begin{subequations}\label{eq:lin_ocp_admm}
    \begin{align}
        \mathbf{r}^{k+1} &\coloneqq \argmin_{\mathbf{r}} \ C_x(\mathbf{r}) + \frac{\rho}{2}\Vert \mathbf{x}^k - \mathbf{r} + \mathbf{v}^k \Vert_2^2 \nonumber \\
        &\quad \text{s.t. } \mathbf{r} \in \mathcal{R}^N \label{eq:lin-ref-layer} \\
        \left(\mathbf{x}^{k+1}, \mathbf{u}^{k+1}\right) &\coloneqq \argmin_{\mathbf{x, u}} \ \frac{\rho}{2}\sum_{t=0}^{N}\Vert x_t - r_t^{k+1} + v_t^k \Vert_2^2 + \sum_{t=0}^{N-1} u_t^T R_t u_t  \nonumber \\
        & \quad \text{s.t. } x_{t+1} = A_t x_t + B_t u_t, \ t=0, \dots, N-1  \nonumber\\
        & \quad \quad \ \, x_0 = \xi \label{eq:lin-track-layer}\\
        \mathbf{v}^{k+1} &\coloneqq \mathbf{v}^k + \mathbf{x}^{k+1} - \mathbf{r}^{k+1}. \label{eq:lin-dual-update}
\end{align}
\end{subequations}

We recognize that the feedback control layer update problem~\eqref{eq:lin-track-layer} is an LQR reference tracking problem,
with the reference trajectory $\mathbf{r}^{k+1}$ which can be solved via dynamic programming. 
We first expand the square to isolate the tracking error term $x_t-r_t^{k+1}$ to obtain the following OCP
\begin{equation}\label{eq:lin_ocp}
    \begin{aligned}
    \min_{\mathbf{x, u}} & \ \sum_{t=0}^{N-1} \frac{\rho}{2}\Vert x_t - r_t^{k+1}\Vert_2^2 + \rho (x_t - r_t^{k+1})^T v_t^k + u_t^T R_t u_t + \frac{\rho}{2}\Vert x_N - r_N^{k+1} \Vert_2^2 + \rho (x_N -r_N)^T v_N^k \\
   \text{s.t.} & \quad \ x_{t+1} = A_t x_t + B_t u_t, \ t=0,\dots,N-1 \\
    & \quad \ x_0 = \xi
    \end{aligned}
\end{equation}

Set $e_t:= x_t - r_t^{k+1}$, $\mu_t := (r_t^{k+1}, r_{t+1}^{k+1},\dots,r_N^{k+1}, 0,\dots,0)$, and $z_t = (e_t, \mu_t)$, and define matrices $F$ and $G$ such that $Fz_t = e_t$ and $Gz_t = \mu_t$.  Then setting $\bar Q_t :=(\rho/2) F^T F$, and $q_t:=(\rho/2)F^Tv_t^k$, we can rewrite problem~\eqref{eq:lin_ocp} as
\begin{equation}\label{eq:opt_ctrl2}
\begin{array}{rl}
         \min_{z_t,u_t} & \sum_{t=0}^{N-1} [z_t^T \bar Q_t z_t + 2q_t^T z_t + u_t^T R_t u_t] + z_N^T \bar Q_N z_N \\
     \text{s.t.} & z_{t+1} = \bar A_tz_t + \bar B_tu_t,\ t=0,\dots,N-1,
\end{array}
\end{equation}
for suitably defined matrices $\bar{A}_t$, $\bar{B}_t$, and $z_0$.  This is a finite horizon LQR optimal control problem with quadratic and affine stage-wise cost terms, which can be solved by dynamic programming. 

We consider a cost-to-go function of the form
\begin{equation}\label{eq:cost-to-go}
    V_t(z_t) = z_t^T P_t z_t + 2p_t^T z_t + c_t,
\end{equation}
where $P_t \in \mathbb{R}^{n(N+2) \times n(N+2)}$ , $p_t \in \mathbb{R}^{n(N+2)}$, and $c_t \in \mathbb{R}$. The terminal cost for the augmented system is obtained by setting $P_N = \Bar{Q}_N=(\rho/2) F^T F$, $p_N = q_N = (\rho/2)F^Tv_N^k$, and $c_N = 0$. We solve the Hamilton-Jacobi equation
$$
V_{t}(z_t) = \min_v z_t^T \Bar{Q}_t z_t + 2q_t ^T z_t + v^T R_t v + V_{t+1}(z_{t+1}),
$$
which has a minimizer given by
\begin{align*}
    v^\star &= -(R_t + \Bar{B}_t^T P_{t+1} \Bar{B}_t)^{-1}(\Bar{B}_t^TP_{t+1}\Bar{A}_tz_t + \Bar{B}_t^Tp_{t+1}) \\
    & =: -K_tz_t - \nu_t,
\end{align*}
where we define $K_t := -(R_t + \Bar{B}_t^TP_{t+1}\Bar{B}_t)^{-1}\Bar{B}_t^TP_{t+1}\Bar{A}_t$ and $\nu_t := -(R_t + \Bar{B}_t^T P_{t+1}\Bar{B}_t)^{-1} \Bar{B}_t^Tp_{t+1})$.  

Plugging $v^\star$ into the cost-to-go function and simplifying further, we observe that the recursions for the matrices  
\begin{equation}\label{eq:riccati}
    P_t = \Bar{Q}_t + K_t^T R_t K_t + (\Bar{A}_t-\Bar{B}_tK_t)^T P_{t+1} (\Bar{A}_t-\Bar{B}_tK_{t+1})
\end{equation} 
follow the usual discrete Algebraic Riccati recursion, and
\begin{align}\label{eq:lin-const-cost}
p_t &= q_t + K_{t+1}^TR_t \nu_{t+1} + (\bar A_t-\bar B_tK_{t+1})^T(p_{t+1}-P_{t+1}\bar B_t\nu_{t+1}) \\
c_t & = -p_{t+1}^T\bar B_t(R_t + \bar B_t^T P_{t+1} \bar B_t)^{-1}\bar B_t^T p_{t+1}.
\end{align}

The optimal control action at time $t$ is then specified by $u_t = -K_tz_t - \nu_t$, which is further decomposed as
\begin{align}\label{eq:ctrl-law}
    u_t &=  -K_t z_t - \nu_t\nonumber \\
    & = -K_tF^TFz_t - K_tG^TGz_t - \nu_t \nonumber\\
    & =: -K^{fb}_te_t - K^{ff}_t\mu_t - \nu_t
\end{align} 
This decomposition highlights that the optimal control action $u_t = -K_t z_t - v_t$ is naturally composed of feedforward and feedback terms that drive the system to and stabilize it around the reference trajectory $\vec r^{k+1}$:
\begin{enumerate}[label=\alph*)]
    \item \textbf{Feedforward term $K^{ff}_t\mu_t$:} this term applies control actions to drive the system towards the desired reference trajectory, as encoded in the look-ahead state $\mu_t =(r^{k+1}_t, \dots, r^{k+1}_N,0,\dots,0)$.
    \item \textbf{Feedback term $K^{fb}_te_t$:} this term stabilizes the system around the nominal trajectory by applying a feedback term based on the error $e_t=(x_t-r_t^{k+1})$.
    \item \textbf{Coordination term $\nu_t$:} The correction term $\nu_t$, which can be seen to be a linear function of the dual variable $\vec v^k$, coordinates the feedback layer behavior with that of planning layer, ensuring convergence to zero tracking error (i.e., that $\vec x=\vec r^{k+1}$) as $k\to\infty$.
\end{enumerate}
We emphasize that the 2DOF structure of the controller was not imposed a priori, and rather naturally emerged from the ADMM algorithm applied to solving OCP~\eqref{eq:redundant-ref-var}.  Further, in contrast to prior work~\cite{matni2016theory} that relied on relaxing the original OCP, the 2DOF layered controller obtained here is optimal.

\textbf{Convergence}: If $C_x$ is a closed, proper, and a convex function and $\mathcal R^N$ is a convex set, Assumption~\ref{a1} is satisfied.  Further, if the linear OCP satisfies strong duality, e.g., if Slater's condition holds, then Assumption~\ref{a2} is satisfied.  This is true, if for example, the state constraint $\mathcal R$ is a polytope, or if it contains the origin in its interior.  It therefore follows by Theorem~\ref{thm:admm} that residual, objective, and dual variable convergence are guaranteed.

\subsection{Stochastic linear system with process noise}

We extend the analysis of the previous section to stochastic linear systems of the form 
\begin{equation}\label{eq:slti}
    x_{t+1} = A_tx_t + B_tu_t + H_tw_t,
\end{equation} 
where $w_t \sim \mathcal{N}(0, I), \forall t$ are i.i.d. zero mean Gaussian with identity covariance. We consider the stochastic linear optimal control problem:
\begin{equation}\label{prob:socp}
    \begin{array}{rl}
         \mathrm{minimize}_{\mathbf{x, u}} & \mathbb{E}_w \left[ \mathcal{C}_x(\mathbf{x}) + \sum_{t=0}^{N-1}u_t^T R u_t \right] \\
         \text{s.t.}& x_{t+1} = A_tx_t + B_tu_t + H_tw_t,\\
         & \mathbb E_w x_t \in \mathcal{R}, \ t=0, \dots, N-1,\\
         & x_0 = \xi.
    \end{array}
\end{equation}
We note that enforcing the constraint in expectation, i.e., $\mathbb E_w x_t \in \mathcal{R}$, could be replaced with suitable chance constraints or moment constraints, but we consider this form of stochastic OCP for simplicity.  Prior to applying the approach of the previous section, we recall that due to linear superposition, the evolution of the stochastic dynamics~\eqref{eq:slti} can be decomposed into deterministic and zero-mean stochastic components, i.e., if we write
\begin{equation}\label{eq:slti_decomp}
\begin{array}{rcl}
     x^d_{t+1} &=& A_tx^d_t + B_tu^d_t, \ x^d_0 = \xi, \\
     x^s_{t+1} &=& A_tx^s_t + B_tu^s_t + H_tw_t, \ x^s_0 = 0,
\end{array}
\end{equation}
then $x_t = x^d_t + x^s_t$ and $u_t = u^d_t + u^s_t$, $(x^d_t,u^d_t)$ are deterministic, and $\mathbb E_w x^s_t = 0$, $\mathbb E_w u^s_t = 0$.

We now apply the approach of the previous section, but introduce redundant reference variables to track only the deterministic component of the dynamics, i.e., we consider the equivalent stochastic OCP:
\begin{equation}\label{prob:red_socp}
    \begin{array}{rl}
         \underset{\vec x^d, \vec  x^s, \vec u^d,\vec u^s,\vec r}{\mathrm{minimize}} & \mathbb{E}_w \left[ \mathcal{C}_x(\mathbf{r} + \vec x^s) + \sum_{t=0}^{N-1}(u^d_t+u^s_t)^T R (u^d_t+u^s_t) \right] \\
         \text{s.t.}& x^d_{t+1} = A_tx^d_t + B_tu^d_t, \ x^d_0 = \xi, \\
                    & x^s_{t+1} = A_tx^s_t + B_tu^s_t + H_tw_t, \ x^s_0 = 0,\\
         & r_t \in \mathcal{R}, \ t=0, \dots, N-1,\\
         & \vec r = \vec x^d.
    \end{array}
\end{equation}

In general, the resulting ADMM updates do not have closed-form expressions, although they are convex and can be approximately solved using stochastic gradient methods.  In order to obtain closed-form expressions, we assume that the state utility function is a convex quadratic, i.e., that $\mathcal C_x(\vec x)= \sum_{t=0}^{N-1} x_t^T C_x x_t + c_x^Tx_t$ for $C_x$ a positive semidefinite matrix, and $c_x$ a vector.\footnote{General costs can also be approximated by their 2nd order Taylor series expansion.} 
 We note that in this case, $\mathbb E_w C_x(\vec x) = \mathbb E_w C_x(\vec x^d + \vec x^s) = C_x(\vec x^d) + \sum_{t=0}^{N-1}\mathrm{Tr}C_x \mathbb E_w x^s_t(x^s_t)^T$, i.e., just as the dynamics do, the deterministic and stochastic components of the cost decouple since $\mathbb E_w x^s_t = 0$.

Applying ADMM to the deterministic component of the optimal control problem yields identical iterates to those found in equation~\eqref{eq:lin_ocp_admm}, and the stochastic component reduces to a standard stochastic LQR problem with cost matrices $(C_x, R)$.  The resulting solution thus inherits the 2DOF layered architecture of the deterministic setting, with $u^d_t$ having the feedforward/feedback structure defined in equation~\eqref{eq:ctrl-law}, and $u^s_t=-K_{LQR}x^s_t$, for $K_{LQR}$ the standard LQR controller defined by the solution to the discrete Algebraic Riccati equation defined in terms of cost matrices $(C_x, R)$ and dynamics $(A,B)$.  Thus, by appropriately applying ADMM to solve stochastic OCP~\eqref{prob:socp}, we show that for quadratic state and control costs, a 2DOF layered control architecture with certainty equivalent trajectory generation and feedback control is optimal.

\section{Extensions}\label{sec:non-lin}

\subsection{Low-order reference trajectories}
In the above, we enforced that $\vec r = \vec x$, i.e., we introduced a reference trajectory of the same dimension as the original state.  In practice, planning is often done using a lower-order reference trajectory such that $r_t = C x_t$ for $C \in \R^{q \times n}$ with $q < n$.  For example, in robotics applications, the state $x = (q, \dot q)$ is composed of generalized coordinates and velocities, it is common to plan only in $r = q$ coordinates.  This is trivially incorporated in the above framework by suitably modifying the redundant equality constraint to enforce $\vec r = \vec C \vec x$ and subsequently applying ADMM.

\subsection{Input constraints}
In the above, we did not consider constraints on the control input of the form $u_t\in\mathcal{U}$, for $\mathcal U$ a convex control input constraint set. We note however that by similarly introducing a redundant control action variable constrained to satisfy $\vec a = \vec u$, and enforcing that $a_t \in \mathcal U$ in the trajectory generation layer~\eqref{eq:ref-layer} problem, now over decision variables $(\vec r, \vec a)$, will ensure input constraint satisfaction.

\subsection{Layered control architectures for nonlinear systems}
We now revisit the general nonlinear OCP~\eqref{prob:ocp} and corresponding ADMM iterate updates~\eqref{eq:admm_ocp}.  If the cost functions $\mathcal{C}_x$ and $\mathcal{C}_u$, as well as the constraint set $\mathcal{R}$, are convex, then the only nonconvex component of the problem is the \emph{unconstrained} nonlinear optimal control problem found in the feedback control layer~\eqref{eq:track-layer}.  By isolating the nonconvexity of the problem to this update step, we can leverage existing techniques from nonlinear optimization and optimal control to approximately solve this update step by applying e.g., iLQR~\cite{todorov2005generalized}, which is guaranteed to rapidly converge to a locally optimal solution under fairly benign assumptions~\cite{liao1991convergence}.In the next section, we demonstrate the usefulness of this decoupling of unconstrained nonlinear optimal control and constrained planning by empirically demonstrating that our ADMM-based LCA converges to better solutions more reliably than vanilla iLQR.  We also show that the modularity of the approach enables more complex constraints, such as obstacle avoidance encoded via integer programming, to be seamlessly integrated into the trajectory generation layer subproblem~\eqref{eq:ref-layer}.

\section{Numerical Examples}\label{sec:sim}

In this section, we present experiments\footnote{All code needed to reproduce these experiments can be found at \url{https://github.com/Nusha97/Layered-control-architectures-for-Robotics/tree/main/dual-ascent}} on a $2-$D linear system in both deterministic and stochastic settings and three nonlinear systems to evaluate our proposed methods. In all of the following, we use the $\rho$-update rule described in~\cite[\S3.4.1, equation (3.13)]{boyd2011distributed} to improve convergence of the ADMM algorithm, and use \texttt{cvxpy}~\cite{diamond2016cvxpy} for solving the convex $\vec r$-update problems. For more details on experiment design, refer to Appendix~\ref{sec:appendix}.

\subsection{2-D linear system}

Consider the discrete-time linear time-invariant system:
\begin{equation}\label{eq:linsys}
x_{t+1}=\begin{bmatrix}
        1 & 1 \\
        0 & 1
    \end{bmatrix}x_t+\begin{bmatrix}
        1 & 0 \\
        0 & 1
    \end{bmatrix}u_t +\begin{bmatrix}
        0.1 & 0 \\
        0 & 0.1
    \end{bmatrix}w_t
\end{equation}
where $x_t, u_t, w_t \in \mathbb{R}^2$, and the initial state $x_0 = (0, 0)$. We seek to design a control law such that the system tracks a circular reference trajectory specified by $s_{1, t} = 2 \cos{\omega t}, s_{2, t} = 2 \sin{\omega t}$ for $\omega = 0.5$, starting from the initial state $x_0$. To this end, we specify the utility cost at the trajectory generation layer as $C_x(\mathbf{r}) = \sum_{t=0}^N \Vert r_t - s_t \Vert_2^2$ over $N=20$ time steps. 

\begin{figure}
    \begin{subfigure}{0.5\textwidth}
    \includegraphics[height=2.4in]{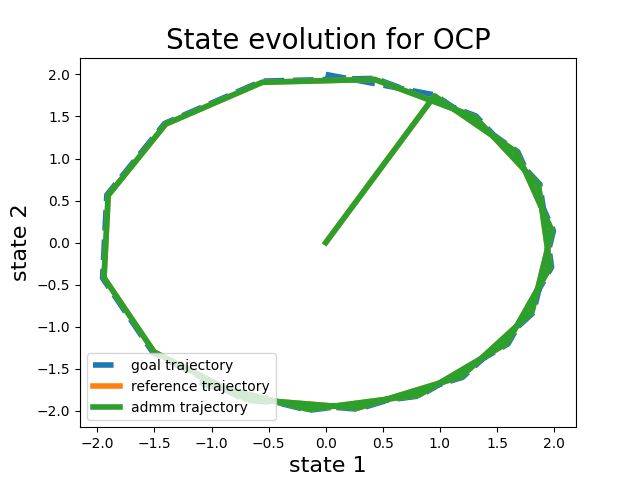}
        \caption{}
    \label{fig:lin_admm}
    \end{subfigure}%
    \begin{subfigure}{0.5\textwidth}
        \includegraphics[height=2.4in]{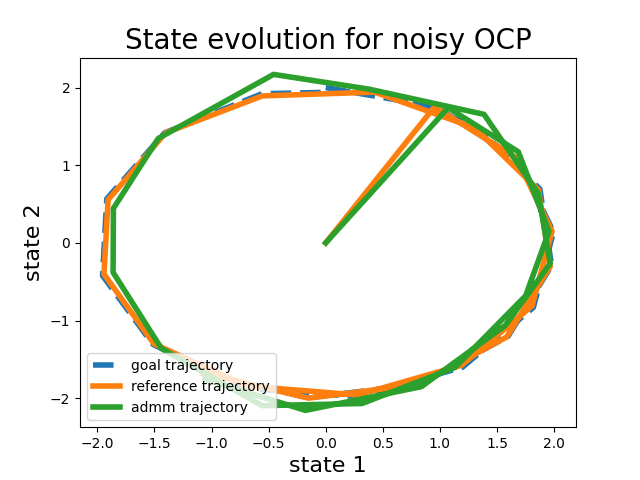}
        \caption{}
        \label{fig:admm_noise}
    \end{subfigure}
    \caption{Plot shows results from running our proposed approach for the linear system~\eqref{eq:linsys} in the absence of disturbances in (a) and in the presence of disturbances in (b). 
    The blue trajectory denotes the reference states $(s_{1, t}, s_{2, t})$, and the orange and green trajectories are respectively the final reference $\mathbf{r}$ and state $\mathbf{x}$ iterates from using the obtained feedback law.}
\end{figure}

\textbf{Results}: As shown in Figure~\ref{fig:lin_admm}, our approach recovers the optimal solution in the deterministic setting when $w_t = 0, \forall t$. In the disturbance setting, our approach as shown in Figure~\ref{fig:admm_noise} recovers the deterministic solution and provides a feedback controller with feedforward and feedback structure. The feedback component from the controller stabilizes the system in the presence of disturbances. We emphasize that in the linear setting, the behavior of these controllers is expected, and hence no exhaustive evaluations or comparisons are necessary.  It is rather the 2DOF structure that emerges as a property of the solution that is of interest here.

\subsection{Nonlinear systems}\label{exp:non-lin}

To evaluate our proposed methods for nonlinear systems, we consider three case studies. The first task is the stabilization of pole dynamics on a moving cart, the second is the navigation of a non-holonomic car-like robot to a goal in the presence of corridor, input and nonconvex obstacle contraints and the last task is to navigate a quadrotor to a goal in 3D.  In the following, we use the \texttt{trajax}~\cite{frostigtrajax} implementation of iLQR. We reserve the term ``convergence" if there is primal and dual feasibility from running our ADMM-based nonlinear LCA. We reserve the term ``success rate" to measure the percentage of trials in which the terminal state from the executed trajectory reaches within a ball of radius $0.5$ from the goal. Refer to Appendix~\ref{sec:appendix} for additional details on experiment design. 

\textbf{Cartpole}: We consider the stabilization task of the pole on a cart and compare the performance of our approach with iLQR~\cite{todorov2005generalized}. The continuous dynamics of the cartpole can be written as:
\begin{equation}
    \begin{aligned}
        H(q)\ddot{q} + C(Q, \dot{q}) \dot{q} + G(q) = F
    \end{aligned}
\end{equation}
where $q = (x, \theta)$ and 
$$H = \begin{bmatrix}
    m_c + m_p & m_p l \cos{\theta} \\
    m_p l \cos{\theta} & m_p l^2
\end{bmatrix}, \ C = \begin{bmatrix}
    0 & -m_p l \dot{\theta} \sin{\theta} \\
    0 & 0
\end{bmatrix},\  G = \begin{bmatrix}
    0 \\
    m_p g l \sin{\theta}
\end{bmatrix}.$$ The system state is given by $(q, \dot{q})\in\R^4$ and the control input is the force applied to the cart in the horizontal direction. We apply Euler discretization to the continuous time system with sampling time $dt = 0.1$, and use the discrete-time dynamics for the rest of the evaluation. 

We apply the approach proposed in Section~\ref{sec:non-lin} using iLQR to solve the feedback control layer problem~\eqref{eq:track-layer} until the ADMM algorithm has converged. For the iLQR step in our ADMM algorithm, we set the maximum number of iterations to $10$, resulting in only approximate (locally optimal) solutions at each iterate update. To compare the performance of our approach, we run iLQR~\cite{todorov2005generalized} to stabilize the pole around the equilibrium point at $(0, \pi, 0, 0)$ with the maximum number of iterations set to $200$.

\begin{figure}
    \begin{subfigure}{0.5\textwidth}
    \includegraphics[height=2.4in]{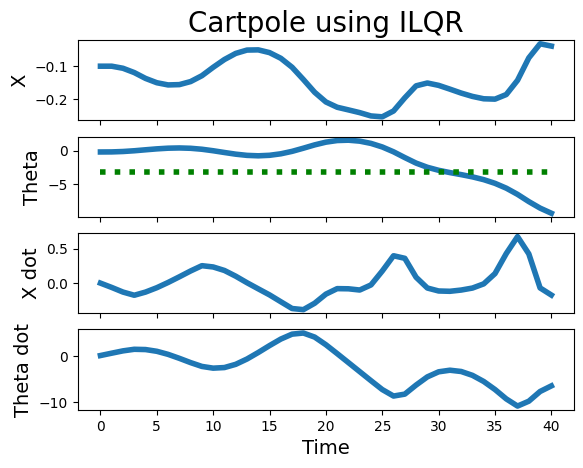}
        \caption{}
    \label{fig:cart_ilqr}
    \end{subfigure}%
    \begin{subfigure}{0.5\textwidth}
        \includegraphics[height=2.4in]{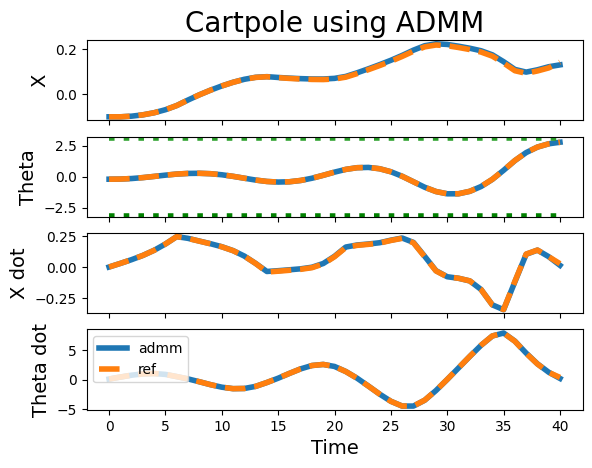}
        \caption{}
        \label{fig:cart_admm}
    \end{subfigure}
    \caption{Representative system traces of running iLQR (a) and ADMM (b) to balance a pole on a cart from the same random initial condition.  Here, iLQR struggles to converge to a stabilizing solution, whereas ADMM is able to successfully achieve balancing in as few as $15$ iterations.}
\end{figure}

\textbf{Results}: We sampled $20$ random initial conditions from a standard uniform distribution and tested our approach against iLQR over a horizon of $40$ time steps and report the results in Table~\ref{tab:cp}. We note that our approach converged to the equilibrium point for every initial condition while iLQR successfully reaches within the goal radius only for $2$ out of $20$ trials.  Qualitatively representative traces of an initial condition for which iLQR fails but our approach succeeds are found in Figs.~\ref{fig:cart_ilqr} and~\ref{fig:cart_admm}. We compute the total number of iterations for our approach as $\sum_{it=1}^K (1 + i_{it})$ 
where $K$ is the number of ADMM outer loops, and $i_{it}$ is the number of iLQR iterations per update step where $i_{it} \leq 10$. We also count every $\vec r$-update step, and thus we are comparing the number of convex optimization oracle calls needed by each algorithm.  As can be observed, the ADMM algorithm appears to demonstrate more favorable performance properties, but as expected, requires more iterations.  We note however that we use a stopping criterion of primal residual error $\|\vec x - \vec r\|_2^2 \leq 10^{-2}$, but we observed that this can be substantially relaxed while still yielding acceptable solutions: we leave optimizing the algorithm for computational efficiency for future work.

\begin{table}[h]
    \centering
    \begin{tabular}{c|c|c}
    \hline
        & No. of iterations &  Success rate (\%) \\
        \hline
        \textit{Ours} & $428.4 \pm 328.7$ & $100$ \\
        \hline
        \textit{iLQR} & $52.5 \pm 41.4$ & $10$ \\ 
        \hline
    \end{tabular}
    \caption{\textmd{Comparison of our approach against iLQR for $20$ randomly sampled initial conditions on cartpole dynamics.}}
    \label{tab:cp}
\end{table}

\textbf{Unicycle}: Next, we consider the navigation task of the continuous time unicycle dynamics given by
$$\begin{bmatrix}
    \dot x_1 \\ \dot x_2 \\ \dot \theta
\end{bmatrix} = \left[\begin{array}{cc}
        \cos{\theta} & 0 \\
        \sin{\theta} & 0 \\
        0 & 1
    \end{array}\right] \left[\begin{array}{c}
         v \\
         \omega 
    \end{array}\right]
$$ where $(x_1,x_2)\in\R^2$ are the system's Cartesian coordinates,  $\theta$ is the heading angle, and $v$, $\omega$ are the instantaneous linear and angular velocities, respectively. We apply Euler discretization to the continuous time dynamical system for the rest of the evaluation as done previously for the cartpole. In addition to dynamics constraints, we include corridor constraints as state constraints and linear speed constraints as input constraints in the nonlinear optimal control problem. 

We apply our proposed approach from Section~\ref{sec:non-lin} using iLQR to solve the feedback control layer problem~\eqref{eq:track-layer} and set the maximum number of iterations to $10$. We compare the performance of our approach against iLQR~\cite{todorov2005generalized} solving the global nonlinear optimal control problem to navigate the car-like robot to a fixed goal at $(3, 2)$.

\begin{figure}
    \begin{subfigure}[t]{0.5\textwidth}
        \includegraphics[height=2.4in]{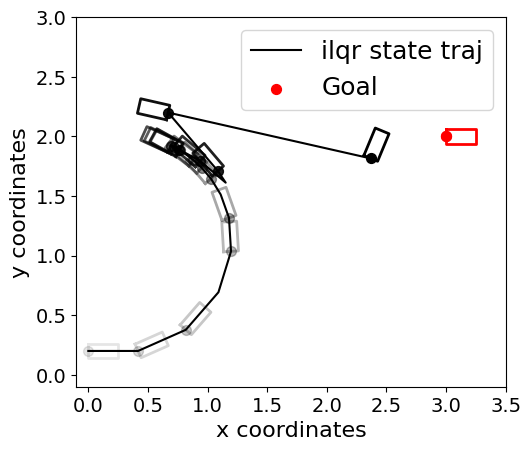}
        \caption{}
    \label{fig:uni_ilqr}
    \end{subfigure}%
    \begin{subfigure}[t]{0.5\textwidth}
        \includegraphics[height=2.4in]{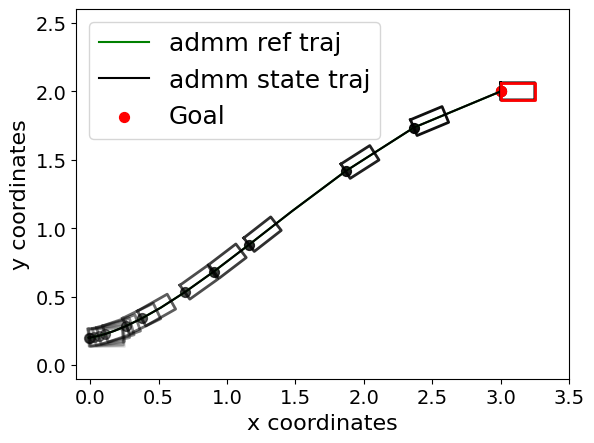}
        \caption{}
        \label{fig:uni_admm}
    \end{subfigure}
    \caption{Plots showing the state trajectories from running iLQR (a) and ADMM (b) for navigating the car-like robot to a goal.}
\end{figure}

\begin{figure}[t]
    \centering
    \includegraphics[height=2.4in, width=0.5\columnwidth]{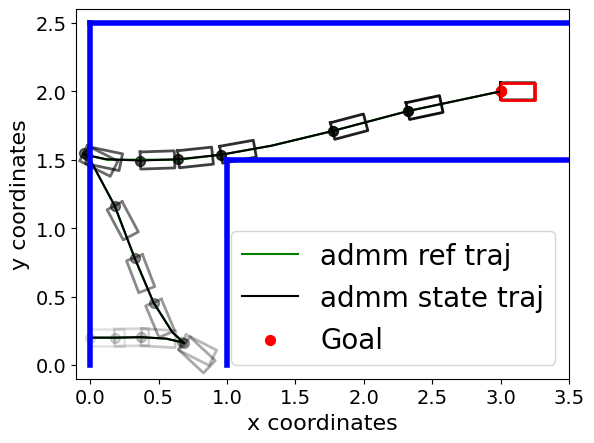}
    \caption{Plot shows the converged reference and state trajectories from running our proposed approach with corridor constraints as shown in blue.}
    \label{fig:admm_corridor}
\end{figure}

\begin{figure}
    \begin{subfigure}[t]{0.5\textwidth}
        \includegraphics[height=2.4in]{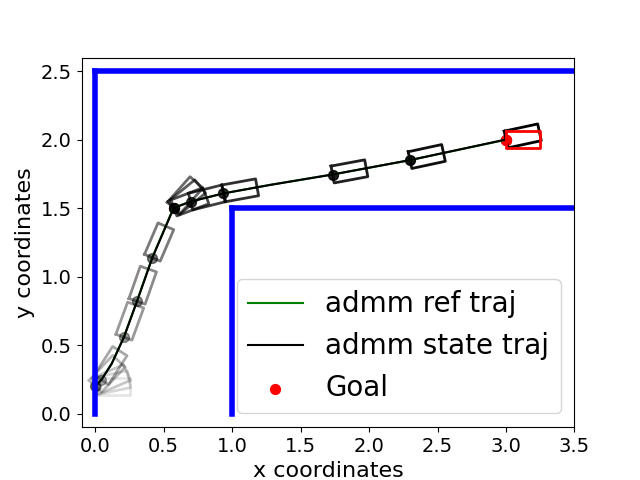}
        \caption{}
    \label{fig:uni_vel_traj}
    \end{subfigure}%
    \begin{subfigure}[t]{0.5\textwidth}
        \includegraphics[height=2.4in]{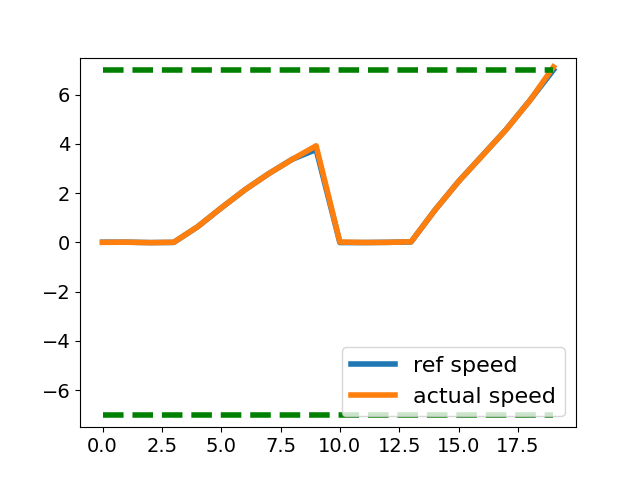}
        \caption{}
        \label{fig:uni_vel_const}
    \end{subfigure}
    \caption{Results from running our proposed approach with a lower order reference in the presence of corridor and velocity constraints. On the left (a), we show the converged reference and state trajectories from navigating the car-like robot to a goal and on the right (b), we plot the linear speed constraints in green dotted lines.}
\end{figure}

\textbf{Results}: We sampled $20$ random initial conditions from a standard normal distribution and tested our approach against iLQR over a horizon of $20$ time steps. As discussed in Section~\ref{sec:non-lin}, we repeat our experiments with lower order reference trajectories planned over $x, y$ positions, and include input constraints on maximum forward and reverse speeds while solving the trajectory generation layer~\eqref{eq:ref-layer}. Additionally, we also evaluated our approach by applying corridor state constraints by switching between affine constraints and also nonconvex obstacle constraints~\ref{sec:appendix} in the planning layer~\eqref{eq:ref-layer}. We report the results in Table~\ref{tab:uni} for the different test cases. We note that our approach, in both cases with and without constraints, successfully reached the goal for every initial condition. In contrast, iLQR successfully reached within the goal radius only $4$ out of $20$ trials without state or input constraints. We show a qualitative comparison of our approach and iLQR in Figs.~\ref{fig:uni_ilqr} and~\ref{fig:uni_admm} respectively where iLQR fails while our approach succeeds in reaching the goal. Next, we show our approach in Figs.~\ref{fig:admm_corridor} and~\ref{fig:uni_vel_traj} finding a feasible path satisfying the non-holonomic behavior of the system, while bringing it to the goal in the presence of tight corridor constraints. In addition, we plot the reference and actual velocities in Fig.~\ref{fig:uni_vel_const} to show the satisfaction of input constraints.

\begin{figure}
    \begin{subfigure}[t]{0.5\textwidth}
        \includegraphics[height=2.4in]{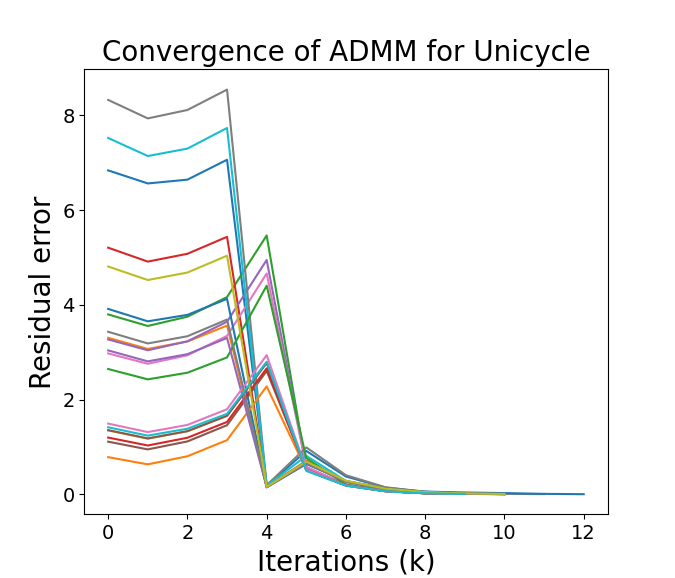}
        \caption{}
    \label{fig:admm_convg}
    \end{subfigure}%
    \begin{subfigure}[t]{0.5\textwidth}
        \includegraphics[height=2.4in]{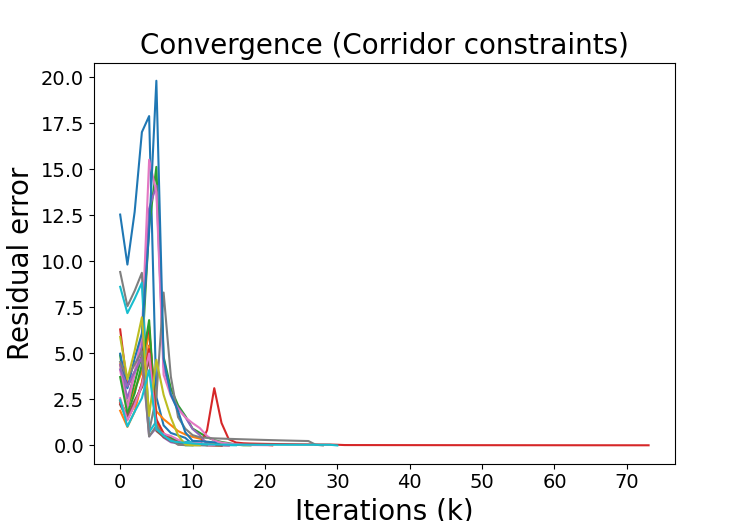}
        \caption{}
        \label{fig:admm_corridor_convg}
    \end{subfigure}
    \caption{We plot the convergence of ADMM on $20$ randomly sampled initial conditions for the unicycle in (a) and for the unicycle dynamics with tight corridor constraints in (b).}
\end{figure}

We do a similar computation of the total number of iterations as explained previously and plot the primal residual error for our approach as shown in Figs.~\ref{fig:admm_convg} and~\ref{fig:admm_corridor_convg}. Using a stopping criterion of primal residual error $\|\vec x - \vec r\|_2^2 \leq 10^{-2}$, we observed that our method converges in as few as $12$ ADMM outer-loop iterations when there are no state constraints and around $26$ iterations on average with corridor constraints. 

\begin{table}[h]
    \centering
    \begin{tabular}{c|c|c}
    \hline
        & No. of iterations &  Success rate (\%) \\
        \hline
        \textit{Ours} & $117.1 \pm 6.3$ & $100$ \\
        \hline
        \textit{Ours (corridor)} & $214.1 \pm 108.3$ & $100$ \\
        \hline
        \textit{Ours (low order)} & $113.8 \pm 5.9$ & $100$ \\
        \hline
        \textit{Ours (low order corr.)} & $152.5 \pm 51.5$ & $100$ \\
        \hline
        \textit{Ours (low order, corr., vel constr)} & $386.7 \pm 74.1$ & $100$ \\
        \hline
        \textit{iLQR} & $6.6 \pm 3.9$ & $20$ \\ 
        \hline
    \end{tabular}
    \caption{\textmd{Comparison of our approach against iLQR for $20$ randomly sampled initial conditions on unicycle dynamics.}}
    \label{tab:uni}
\end{table}

\textbf{Quadrotor}: Lastly, we consider the navigation task of a quadrotor to a goal with control-affine dynamics from~\cite[Ch. 2]{sabatino2015quadrotor} where the quadrotor states are given by its position, velocity, roll, pitch, yaw and angular velocity in the world frame, and control inputs are given by the collective motor thrusts and body moment torques. Making a simplifying assumption that the angular body rates are equal to the angular velocity in the world frame leads to a control-affine system. We evaluate our proposed approach against iLQR~\cite{todorov2005generalized} to navigate the quadrotor to a fixed goal at $(3, 2, 1.5)$.

\begin{figure}
    \begin{subfigure}[t]{0.5\textwidth}
        \includegraphics[height=2.4in]{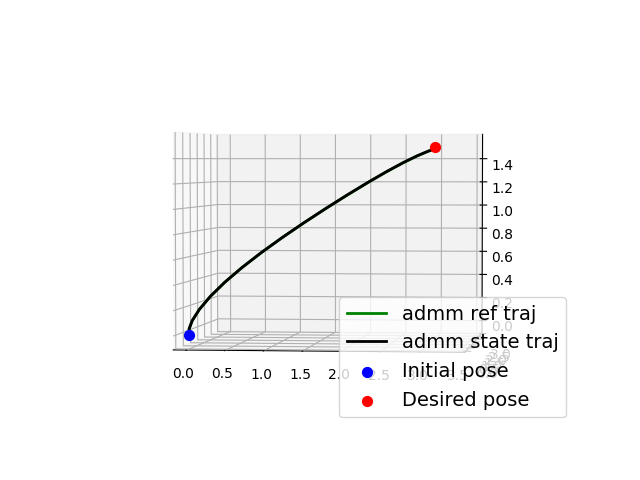}
        \caption{}
    \label{fig:quad_full_order}
    \end{subfigure}%
    \begin{subfigure}[t]{0.5\textwidth}
        \includegraphics[height=2.4in]{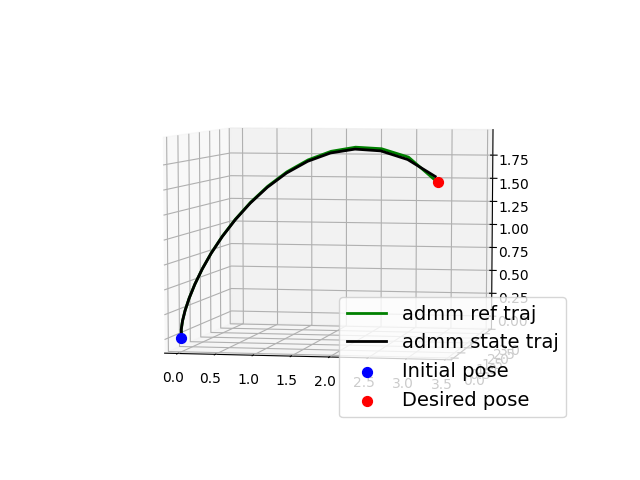}
        \caption{}
        \label{fig:quad_low_order}
    \end{subfigure}
    \caption{On the left (a) we show the plot of the converged 3D reference and state trajectories from running our approach using full order reference, and on the right (b) we show results from using a lower order reference trajectory.}
\end{figure}

\textbf{Results}: We simulate with a horizon of $20$ time steps by sampling $20$ initial conditions from a standard normal distribution. We report the results in Table~\ref{tab:quad} from running our approach on the full order and lower order(consisting of $x, y, z$ positions) reference trajectories. We note that iLQR failed to reach the goal for all initial conditions while our approach found a dynamically feasible trajectory to the goal for both the lower order and full order trajectory planning problems as shown in Figs.~\ref{fig:quad_full_order} and~\ref{fig:quad_low_order}. We leave optimizing our approach for computational efficiency and including state and input constraints from real hardware platforms for future work.

\begin{table}[h]
    \centering
    \begin{tabular}{c|c|c}
    \hline
        & No. of iterations & Success rate (\%) \\
        \hline
        \textit{Ours} & $500.5 \pm 83.6$ & $100$ \\
        \hline
        \textit{Ours (low order)} & $20103.6 \pm 40500$ & $100$ \\
        \hline
        \textit{iLQR} & $4.5 \pm 1$ & $0$ \\
        \hline
    \end{tabular}
    \caption{\textmd{Comparison of our approach against iLQR for $20$ randomly sampled initial conditions on the quadrotor dynamics model.}}
    \label{tab:quad}
\end{table}

\vspace{-2em}

\section{Conclusion}
\label{sec:conclusion}

We showed that by introducing a redundant reference variable to an optimal control problem and subsequently applying ADMM, optimal controllers with a layered structure are obtained.  We instantiated this approach in the context of linear optimal control problems, and recovered a feedforward/feedback-based optimal controller.  In the context of nonlinear optimal control, we empirically demonstrated the benefits of separating trajectory generation from feedback control in terms of both convergence (as compared to vanilla iLQR) and flexibility (by seamlessly incorporating both convex and nonconvex constraints).  Exciting directions of future work include developing convergence guarantees for the proposed nonlinear control scheme by making connections to the nonconvex ADMM literature~\cite{wang2019global}, as well as considering alternative planning problems, such as those based on semantic specifications~\cite{kress2009temporal}.

\appendix
\section{Implementation Details}
\label{sec:appendix}

\subsection{Linear and nonlinear system control design}
For our experiments on linear systems, we obtain a controller by solving an LQR problem in the feedback control layer~\eqref{eq:lin-track-layer} as shown in~\eqref{eq:opt_ctrl2} by specifying the control weight matrix $R_t = 0.001I$ over $t=0, 1, \cdots, N-1$. The state weight matrices $\bar Q_t :=(\rho/2) F^T F$ and $q_t:=(\rho/2)F^Tv_t^k$ use $F \in \mathbb{R}^{n \times (N+2)}$, obtained by horizontally stacking the identity matrix with zeros. 

For our experiments in Section~\ref{exp:non-lin}, we choose the penalty $Q_{r_t}$ on reference states in the trajectory generation layer~\eqref{eq:track-layer} as $Q_{r_t} = 0.1I, \forall t=\{0, 1, \cdots, N-1\}$ and the penalty of deviation of the terminal state from the goal as $Q_{r_N} = 1000I$. To solve the feedback control layer~\eqref{eq:track-layer} using iLQR, we choose control cost matrices as $R_t = 0.01I, \forall t$. We initialize $\rho=25$. We use the same design parameters for the stage, terminal and input costs for iLQR.

\subsection{Low order reference trajectory design}

As mentioned in Section~\ref{sec:non-lin}, reference planning typically uses a lower dimensional state representation which we term as a lower order reference trajectory. We accommodate the lower order reference trajectory by solving the following optimization sub-problems:
\begin{subequations}
    \begin{align}
        \mathbf{r}^{k+1} &\coloneqq \argmin_{\mathbf{r}} \ C_x(\mathbf{r}) + \frac{\rho}{2}\Vert \mathbf{C} \mathbf{x}^k - \mathbf{r} + \mathbf{v}^k \Vert_2^2 \nonumber \\
        &\quad \text{s.t. } \mathbf{r} \in \mathcal{R}^N \\
        \left(\mathbf{x}^{k+1}, \mathbf{u}^{k+1}\right) &\coloneqq \argmin_{\mathbf{x, u}} \ \frac{\rho}{2}\sum_{t=0}^{N}\Vert C x_t - r_t^{k+1} + v_t^k \Vert_2^2 + \sum_{t=0}^{N-1} u_t^T R_t u_t  \nonumber \\
        & \quad \text{s.t. } x_{t+1} = f(x_t, u_t), \ t=0, \dots, N-1  \nonumber\\
        & \quad \quad \ \, x_0 = \xi \\
        \mathbf{v}^{k+1} &\coloneqq \mathbf{v}^k + \mathbf{Cx}^{k+1} - \mathbf{r}^{k+1}. 
\end{align}
\end{subequations}

We note here that for the unicycle model experiments, $r_t \in \mathbb{R}^2$ denotes the $x, y$ positions on the 2D plane and for the quadrotor model experiments, $r_t \in \mathbb{R}^3$ denotes the $x, y, z$ positions in 3D. 

\subsection{Corridor and input constraint design}

For the design of input constraints as described in Section~\ref{sec:non-lin}, we solve the following optimization sub-problems:
\begin{subequations}
    \begin{align}
        \mathbf{(r^{k+1}, a^{k+1})} &\coloneqq \argmin_{\mathbf{r, a}} \ C_x(\mathbf{r}) + \frac{\rho}{2}\Vert \mathbf{x}^k - \mathbf{r} + \mathbf{v_r}^k \Vert_2^2 + \frac{\rho}{2}\Vert \mathbf{u}^k - \mathbf{a} + \mathbf{v_a}^k \Vert_2^2 \nonumber \\
        &\quad \text{s.t. } \mathbf{r} \in \mathcal{R}^N, \quad a_t \in \mathcal{U} \\
        \left(\mathbf{x}^{k+1}, \mathbf{u}^{k+1}\right) &\coloneqq \argmin_{\mathbf{x, u}} \ \frac{\rho}{2}\sum_{t=0}^{N}\Vert x_t - r_t^{k+1} + v_{r_t}^k \Vert_2^2 + \sum_{t=0}^{N-1} u_t^T R_t u_t + \frac{\rho}{2} \Vert u_t - a_t^{k+1} + v_{a_t}^k \Vert_2^2 \nonumber \\
        & \quad \text{s.t. } x_{t+1} = f(x_t, u_t), \ t=0, \dots, N-1  \nonumber\\
        & \quad \quad \ \, x_0 = \xi \\
        \mathbf{v_r}^{k+1} &\coloneqq \mathbf{v_r}^k + \mathbf{x}^{k+1} - \mathbf{r}^{k+1} \nonumber \\
        \mathbf{v_a}^{k+1} &\coloneqq \mathbf{v_a}^k + \mathbf{u}^{k+1} - \mathbf{a}^{k+1}. 
\end{align}
\end{subequations}

In our unicycle experiments, we limit the maximum linear speeds to be $\pm 7$ m/s in the forward and reverse directions. We choose the constraint sets $\mathcal{U}$ to be linear of the form $|a_t| \leq 7, \forall t$. To enforce corridor constraints as shown in Figs.~\ref{fig:admm_corridor} and~\ref{fig:uni_vel_traj}, we switch between different affine constraints based on the corridor that the car-like robot needs to stay within. For the first half of the time horizon, we enforce constraints on $\mathbf{r}$ to be within the left of the boundary given by $x = 1$ and to the right of $x = 0$ leaving the y-coordinates of $\mathbf{r}$ unconstrained. For the second half of the time horizon, we enforce constraints on $\mathbf{r}$ to be above $y=1.5$ and below $y=2.5$ leaving the x-coordinates of $\mathbf{r}$ unconstrained. Essentially, the corridor state constraints are $r_{1, t} \leq 1, r_{1, t} \geq 0, \forall t \in \{0, \cdots, \frac{T}{2}\}$ and $r_{2, t} \leq 2.5, r_{2, t} \geq 1.5, \forall t \in \{\frac{T}{2}+1, \cdots, T\}$. 

On running our approach, we observe that the solver of the trajectory planning layer returns ``infeasible" if the provided input and state constraints cannot be satisfied while simultaneously finding a path to the goal within the given time horizon. We leave exploring the dependence of constraints and planning horizon for future work.

\subsection{Obstacles as integer constraint design}

To avoid obstacles in the environment, we show how integer constraints can be seamlessly integrated into our formulation. We model the obstacle constraints using integer variables for each time step as discussed in~\cite[Sec. 2.3]{richards2002trajectory} for rectangular obstacles. To design an obstacle for our unicycle model, we specify the rectangle using the lower left corner given by $(x_{min}, y_{min}) = (1, 0.5)$ and the upper right corner given by $(x_{max}, y_{max}) = (1.5, 1)$. There are $5$ constraints that we include per time step as shown below:
\begin{align*}
    r_{1, t} &\leq x_{min} + Ma_1 \\
    -r_{1, t} &\leq -x_{max} + Ma_2 \\
    r_{2, t} &\leq y_{min} + Ma_3 \\
    -r_{2, t} &\leq -y_{max} + Ma_4 \\
    &\sum_{k=1}^4 a_k \leq 3
\end{align*}

\begin{figure}[t]
    \centering
    \includegraphics[height=2.4in, width=0.5\columnwidth]{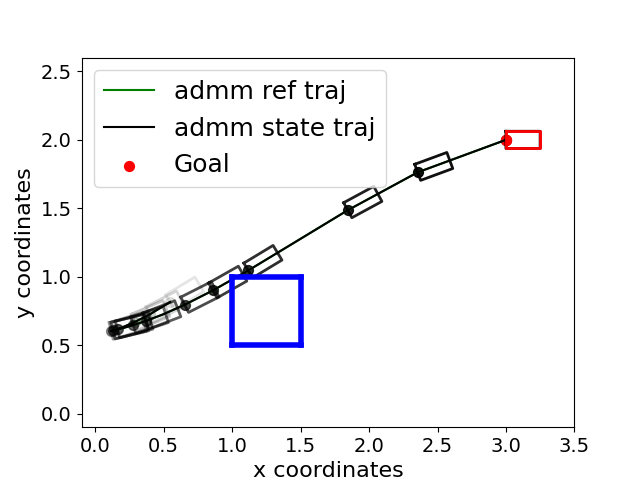}
    \caption{Plot shows the converged reference and state trajectories from running our proposed approach avoiding a rectangular obstacle as shown in blue. We note here that the black dots are the Euler discretized points and integer constraints are imposed on the $x,y$ reference states. Hence, the obtained trajectory is indeed feasible as all the black dots lie outside the blue box.}
    \label{fig:uni_obs}
\end{figure}

We select $M=10^8$ and use \texttt{cvxpy} to solve the trajectory planning layer. We plot the trajectory from running our approach for navigating the car-like robot to a goal from a randomly sampled initial condition in Fig.~\ref{fig:uni_obs}. Although the graphic looks like the car-like robot is moving inside the boundary, we note that the black dots are the Euler discretized points and integer constraints are only imposed on the $x,y$-coordinates of the reference states. Hence, the obtained trajectory is indeed feasible as all the black dots lie outside the blue box. This opens up the possibility of exploring more complex specifications either through discrete constraints or search-based methods. 

\bibliographystyle{abbrvnat}

\bibliography{papers}

\end{document}